\documentclass{elsart}
\usepackage{stmaryrd}
\usepackage{amsfonts}
\usepackage{amsmath}
\usepackage{mathrsfs}
\usepackage{xcolor}
\usepackage[all]{xy}
\usepackage{graphicx}

\allowdisplaybreaks[4]
\newtheorem{tm}{Theorem}[section]
\newtheorem{pn}[tm]{Proposition}
\newtheorem{lm}[tm]{Lemma}

\theoremstyle{definition}
\newtheorem{dn}[tm]{Definition}
\newtheorem{rk}[tm]{Remark}
\newtheorem{ex}[tm]{Example}
\usepackage{enumerate}

\usepackage{bm}

\begin{document}
	
\begin{frontmatter}
\title{The characterization  for the sobriety of    $L$-convex spaces}
\author{Guojun Wu, Wei Yao}
\address{School of Mathematics and Statistics, Nanjing University of Information Science and Technology, Nanjing, 210044, China\\
Applied Mathematics Center of Jiangsu Province, Nanjing University of Information Science and Technology, Nanjing, 210044, China\\
}

\address{}
\date{}

\begin{abstract}
With a commutative integral quantale
$L$ as the truth value table, this study focuses on the  characterizations  of the sobriety of  stratified $L$-convex spaces, as  introduced by Liu and Yue in 2024. It is shown that a stratified sober $L$-convex space $Y$ is a sobrification of a stratified $L$-convex space $X$ if and only if there exists a quasihomeomorphism from $X$ to $Y$;
 a stratified  $L$-convex space is sober if and only if it is a   strictly injective object in the category of  stratified $S_0$ $L$-convex spaces.
 \end{abstract}

\end{frontmatter}

\section{\bf Introduction }
A convex structure, also known as an algebraic closure system, on a set is a family of {\color{black} subsets} closed under arbitrary intersections and directed unions,  including the empty set. A set equipped with a convex structure is called a convex space, or an algebraic closure space. Convex structures can be viewed as  abstractions of traditional convex sets in Euclidean spaces.   Monograph \cite{Vanbook} provides a comprehensive overview of the theory of convex structures.
Convex structures have close connections with various mathematical structures, including
ordered structures  \cite{Erne-stone,Varlet,Van1984,Wu-Guo-Li}, algebraic structures \cite{Liu-Meng-Yue,Algebra,algebra4,Pang-2023},  and topological
structures \cite{topology2,topology1,Yue-Yao-Ho}.

Sobriety is an interesting property in the realm of non-Hausdorff topology \cite{Jean-2013} and domain theory \cite{Gierz}.
It fosters the interplay between ordered and topological structures \cite{XU-mao-form}. The theory of convex spaces also addresses {\color{black}the so-called} sobriety as an important concept.
The notion of sobriety in convex structures can be traced back to the pioneering work of Ern\'{e}
 \cite{Erne-stone},
where he established a Stone-type duality between sober convex spaces and algebraic lattices.
Shen et al. \cite{Shen} explored properties of sober convex spaces similar to those of sober topological spaces and shed light on more  connections between convex spaces and domain theory.
In  \cite{Yao-Zhou}, Yao and Zhou established a categorical isomorphism between   join-semilattices and  sober convex spaces.  Sun and Pang \cite{Sun-Pang-2024}
 proposed an alternative form of sober convex spaces and demonstrated its equivalence to  Ern\'{e}'s sober convex spaces
 \cite{Erne-stone}.
 Xia in \cite{Xia} investigated further properties of sober convex spaces, contributing to the development of pointfree convex geometry research methods \cite{Mar-free-2010}.

With the development of fuzzy set theory and its applications,   {\color{black} fuzzy  convex spaces} have been extensively studied (see \cite{RosaDc,Maruyama,PangShi,Pang2,Su-Li-Intel,Shen1}). Recently,  Liu and Yue  \cite{Liu-Meng-Yue}  extended the concept of algebraic irreducible convex set to the fuzzy setting,   leading to   a fuzzy type  of  sobriety of  $L$-convex spaces.  In  \cite{Wu-Yao}, we extended the notion of  polytopes to the fuzzy setting and introduced an alternative type  of fuzzy sober convex spaces. In \cite{Liu-Meng-Yue}, Liu and Yue   mainly studied the reflectivity and their sobriety in the category of $L$-convex spaces. Different from Liu and Yue's work, we
primarily investigated  $L$-join-semilattice completions  of   $L$-ordered sets via  our notion of sobriety in \cite{Wu-Yao},  establishing more connections  between fuzzy ordered structures and fuzzy convex structures.
 In the classical case, since the notions
of polytopes and algebraic irreducible convex sets are equivalent, one can conclude that these two types of sobriety in classical convex spaces are also equivalent. However, in the fuzzy setting, while Liu and Yue's notion of sobriety clearly implies ours, the reverse implication was left as an unresolved problem in \cite{Wu-Yao}.  {\it {\color{black} Fortunately, in this paper, we provide a counterexample showing that the reverse implication does not hold.}}

 In  topology, it is well known that sober topological spaces
are precisely the strictly injective objects in the category of $T_0$ topological spaces \cite{M.H.Es,Gierz}.  Xia \cite{Xia} provided a corresponding result for sober convex spaces.    In the theory of fuzzy convex structures, it is natural to ask {\it {\color{black} whether sober $L$-convex spaces are exactly the strictly injective objects in the category of $S_0$ $L$-convex spaces}}. This paper will give a positive answer. Moreover, Xia in \cite{Xia} presented a characterization for the sobrification of convex spaces with help of the notion of quasihomeomorphism.  {\it {\color{black} As an application of   fuzzy orders, we also give a characterization for the sobrifications of $L$-convex spaces and  provide an example to demonstrate that this result is helpful in identifying the sobrifications of $L$-convex spaces.}}

In this paper,
we select  a commutative integral quantale $L$ as the truth value table. The  paper is structured  as follows:  In Section 2, we review fundamental concepts and results related to $L$-orders and $L$-convex spaces.  In Section 3, we
revisit essential findings on sober
$L$-convex spaces and provide additional propositions concerning sobriety.
  In Section 4,   we introduce the concepts of quasihomeomorphism and strict embedding within the framework of  $L$-convex spaces, offering characterizations for the sobrification and sobriety of an
$L$-convex space, respectively.

\section{\bf Preliminaries}
We refer to \cite{Quantale} for contents of quantales,   to \cite{MVTop} for notions of fuzzy sets, and to \cite{PartI,Yao-Lu,FanZhang} for contents of fuzzy posets.

Let $L$ be a complete lattice with  a bottom element $0$ and a top element 1 and let $\otimes$ be a binary operation
on $L$ such that $(L,\otimes, 1)$ forms  a commutative monoid.  The pair $(L,\otimes)$ is called a {\it commutative integral quantale}, or a {\em complete residuated lattice},
if the operation $\otimes$ is distributive over joins, i.e.,
$$a\otimes(\bigvee S)=\bigvee_{s\in S} (a\otimes s).$$
For {\color{black} a commutative}  integral quantale $(L,\otimes)$, the operation
$\otimes$ has a corresponding operation  $\rightarrow:L\times L\longrightarrow L$
defined via the  adjoint property:  $$a\otimes b\leq c\ \Longleftrightarrow\ a\leq b\rightarrow c\ (\forall
a,b,c\in L).$$

 In this paper,   unless otherwise specified, the truth value table $L$ is always  assumed to be  a commutative  integral quantale.\\

Every mapping $A:X\longrightarrow L$ is called an $L$-{\it subset} of $X$, denoted by $A\in L^X$. Let $Y\subseteq X$ and {\color{black} $A\in L^X$. Define }$A|_{Y}\in L^Y$ by $A|_{Y}(y)=A(y)$ $(\forall y\in Y)$.
For an element $a\in L$, the notation $a_X$ denotes the constant $L$-subset of $X$ with  the value $a$, i.e., $a_X(x)=a\ (\forall x\in X)$.
For all $a\in L$ and $A\in  L^{X}$, we  write $a\otimes A$  for the $L$-subset given by $(a\otimes A)(x) = a\otimes A(x)$.

For each $a\in X$,  define the characteristic function  $1_{a, X}:X\longrightarrow L$    by
$$
\ \ 1_{a, X}(x)=\left\{\begin{array}{ll}1,& x=a;\\
0,&  x\neq a.
\end{array}\right.
$$
We always write $1_{a}$ instead of {\color{black}$1_{a, X}$}  when there is no confusion.  Readers can identify the domain of a given characteristic function based on the context.
\begin{dn}{(\cite{PartI,FanZhang})}
A mapping $e:P\times P\longrightarrow L$ is called
 an {\it $L$-order} if

{\rm (E1)} $\forall x\in P$, $e(x,x)=1$;

{\rm (E2)} $\forall x,y,z\in P$, $e(x,y)\otimes e(y,z)\leq e(x,z)$;

{\rm (E3)} $\forall x,y\in P$, if $e(x,y)\wedge e(y,x)=1$,
then $x=y$.

\noindent The pair $(P,e)$ is called an $L$-{\em ordered set}. It is customary to write $P$ for the pair $(P, e)$.
\end{dn}
To avoid confusion, we sometimes use $e_P$ to emphasize that the background set is $P$.
A mapping  $f: P\longrightarrow Q$  between two $L$-ordered sets is said to be {\it $L$-order-preserving}
 if for all $x, y\in P$, $e_{P}(x, y)\leq e_{Q}(f(x), f(y))$; $f$  is said to be  {\color{black}an {\it $L$-order-isomorphism}}
 if $f$ is a bijection and $ e_{P}(x, y)= e_{Q}(f(x), f(y))$ for all $x, y\in P$, denoted by $P\cong Q$.

\begin{ex}\label{ex-L-ord}(\cite{PartI,FanZhang})
 Define ${\rm sub}_X :L^X\times L^X\longrightarrow L$ by
 $${\rm sub}_X(A,B)=\bigwedge\limits_{x\in X}A(x)\rightarrow B(x)\ (\forall A, B\in L^X).$$ Then ${\rm sub}_X$ is an $L$-order on $L^X$, which is called the  {\em   inclusion $L$-order} on $L^X$.
If the background set $X$ is clear,  we can  delete the subscript $X$ in ${\rm sub}_X$.
\end{ex}
     A mapping $f:X\longrightarrow Y$ can induce two mappings naturally; those are, $f^\rightarrow:L^X\longrightarrow L^Y$ and $f^\leftarrow:L^Y\longrightarrow L^X$, which are respectively defined by
\vskip 4pt
\centerline{$f^\rightarrow(A)(y)=\bigvee\limits_{f(x)=y}A(x)\ (\forall A\in L^X)$,\   \ \qquad $f^\leftarrow(B)=B\circ f\ (\forall B\in L^Y).$}
\vskip 5pt
\noindent These two mappings are called the {\it Zadeh extensions} of $f$.
\begin{lm} \em(\cite{Yao-Lu}) For each mapping $f: X\longrightarrow Y$,
 {\color{black} the Zadeh extensions} $f^{\rightarrow}$ and $f^{\leftarrow}$are $L$-order-preserving with respect the inclusion $L$-order, and
 $f^{\rightarrow}$ is left adjoint to $f^{\leftarrow}$, denoted   $f^{\rightarrow}\dashv f^{\leftarrow}$; this means that
$$ {\rm sub}_Y(f^{\rightarrow}(A), B)={\rm sub}_X(A, f^{\leftarrow}(B))\ (\forall A\in L^X, B\in L^Y). $$
\end{lm}


Define ${\uparrow}x$ and ${\downarrow}x$ respectively by ${\uparrow}x(y)=e(x,y)$, ${\downarrow}x(y)=e(y,x)\ (\forall x,y\in P)$.
\begin{dn} {(\cite{PartI,FanZhang})}\label{dn-sup-inf} Let $P$ be an $L$-ordered set. An element $x_0$ is called a {\it supremum} of an $L$-subset $A$ of $P$, in symbols $x_0=\sqcup A$, if

{\rm (1)} $\forall x\in P,\ A(x)\leq e(x,x_0)$;

{\rm (2)} $\forall y\in P,\ \bigwedge\limits_{x\in P}A(x)\rightarrow e(x,
y)\leq e(x_0,y)$.
\\An element $x_0$ is called an {\it infimum} of an $L$-subset $A$ of $P$, in symbols $x_0=\sqcap A$,
if

{\rm (3)}  $\forall x\in P,\ A(x)\leq e(x_0,x)$;

{\rm (4)} $\forall y\in P,\ \bigwedge\limits_{x\in P}A(x)\rightarrow
e(y,x)\leq e(y,x_0)$.
\end{dn}

It is easy to verify  that  {\color{black} if the supremum (resp., infimum) of an $L$-subset in an $L$-ordered set exists, then it must be unique.} The following provides an equivalent definition of  supremum and infimum, respectively.

\begin{lm}{\em (\cite{PartI,FanZhang})} Let $P$ be an $L$-ordered set. An element $x_0\in P$ is the  supremum (resp., infimum) of $A\in L^P$, i.e., $x_0=\sqcup A$ (resp., $x_0=\sqcap A$), iff
\vskip 6pt
\centerline{$e(x_0,y)={\rm sub}(A,{\downarrow}y)\ (\forall y\in P)\ (resp., e(y,x_0)={\rm sub}(A,{\uparrow}y)\ (\forall y\in P).$}
\vskip 3pt
\end{lm}

An $L$-ordered set $P$ is said to be {\em complete} if its every $L$-subset has a supremum, or equivalently, has an infimum.
 \begin{ex}{ (\cite{PartI})}
For every $L$-subset $\mathcal{A}$ in the $L$-ordered set $(L^X, {\rm sub})$, the supremum (resp., infimum) of $\mathcal{A}$ exists;
that is
\vskip3pt
\centerline{$\sqcup \mathcal{A}=\bigvee_{B\in L^X}\mathcal{A}(B)\otimes B$ (resp., $\sqcap\mathcal{A}=\bigwedge_{B\in L^X}\mathcal{A}(B)\to B$).}
\vskip3pt
\noindent Hence, the $L$-ordered set $(L^X, {\rm sub})$ is  complete.

\end{ex}

 \begin{dn} Let $P$ and $Q$ be two complete $L$-ordered sets.
The mapping $f:P\longrightarrow Q$  {\color{black} is said to be} {\em supremum-preserving} (resp., {\em infimum-preserving}) if for every $A\in L^{P}$, $f(\sqcup A)=\sqcup f^{\rightarrow}(A)$ (resp., $f(\sqcap A)=\sqcap f^{\rightarrow}(A)$).
\end{dn}
 \begin{ex}{ (\cite[Proposition 5.1]{Yao2016},\cite[Theorem 3.5]{Yao-Lu})}
Let $f:X\longrightarrow Y$ be a mapping between two sets. Then  $f^\rightarrow:(L^X, {\rm sub}_X)\longrightarrow (L^Y, {\rm sub}_Y)$ is supremum-preserving, and $f^\leftarrow:(L^Y, {\rm sub}_Y)\longrightarrow (L^X, {\rm sub}_X)$ is both supremum-preserving and infimum-preserving.
\end{ex}

In the following, we  recall the contents of $L$-convex spaces/algebraic $L$-closure  spaces.

\begin{dn}{ (\cite[Definitions 2.4, 2.5]{Liu-Yue-2024})}\label{dn-con} Let $X$ be a set and $\mathcal{C}\subseteq L^X$. The family $\mathcal{C}$ is called an {\em $L$-convex structure}, or  {\em algebraic $L$-closure structure},
on $X$ if it satisfies the following conditions:

{\rm (C1)} $0_X, 1_X\in \mathcal{C}$;

{\rm (C2)}  $\bigvee_{i\in I}^{\uparrow}C_i\in \mathcal{C}$ for every directed subset  $\{ C_i\mid i\in I\}$  of $\mathcal{C}$;

{\rm (C3)}  $\bigwedge_{j\in J} C_j\in \mathcal{C}$ for every subset $\{ C_j\mid j\in J\}$ of $\mathcal{C}$;

The pair $(X, \mathcal{C})$ is called an {\em $L$-convex space}, or  {\em algebraic $L$-closure space}; each element of $\mathcal{C}$ is called a {\em convex set} of $(X, \mathcal{C})$.

The  $L$-convex space $(X, \mathcal{C})$ is said to be {\em stratified} if  it furthermore satisfies:

{\rm (S)}  $p\to C\in \mathcal{C}$ for all $p\in L$ and $C\in \mathcal{C}$.
\end{dn}
In this paper, every $L$-convex space  is assumed to be stratified, so we omit the the word ``stratified''.
As usual, we often write
$X$ instead of $(X,\mathcal{C})$ for an $L$-convex space  and write $\mathcal{C}(X)$ for the $L$-convex structure of $X$.

\begin{dn}{ (\cite{Liu-Yue-2024})} Let $X$ be an $L$-convex space. Define a mapping $co_X: L^X\longrightarrow L^X$ by
$$ co_X(A)=\bigwedge\{C\in \mathcal{C}(X)\mid A\leq C\}\ (\forall A\in L^X),$$
called the {\em hull operator} of $X$ and $co_X(A)$ is called the {\em hull of $A$}. We  write $co$ instead of $co_X$ when no ambiguity can arise.
\end{dn}

%
%
%
%
%
\begin{dn}{(\cite{PangShi})} Let  $f: X\longrightarrow Y$ {\color{black}be} a  mapping between two $L$-convex spaces.
Then $f$ is said to be
 \begin{itemize}
\item  {\em convexity-preserving}  if for every $D \in \mathcal{C}(Y)$, $f^{ \leftarrow}(D)\in \mathcal{C}(X)$;

\item {\em convex-to-convex}  if for every $C\in \mathcal{C}(X)$, $f^{ \rightarrow}(C)\in \mathcal{C}(Y)$;
\item  {\em homeomorphic} if it is bijective, convexity-preserving and convex-to-convex.
 \end{itemize}

\end{dn}

One easily verifies that a bijection  $f: X\longrightarrow Y$ is a homeomorphism if and only if $f: X\longrightarrow Y$ and its inverse mapping $f^{-1}: Y\longrightarrow X$  are all convexity-preserving.

Given a convexity-preserving mapping $f: X\longrightarrow Y$, we obtain an  assignment from  $\mathcal{C}(Y)$ to $\mathcal{C}(X)$, which sends $D \in \mathcal{C}(Y)$
to $f^{ \leftarrow}(D)\in \mathcal{C}(X)$. When no {\color{black} confusion} can arise, we also use
$$f^{ \leftarrow}:\mathcal{C}(Y)\longrightarrow\mathcal{C}(X)$$ to denote this mapping.  For each $L$-convex space $X$, $(\mathcal{C}(X), {\rm sub}_X)$ is a complete $L$-ordered set. Specifically, for every $L$-subset $\mathcal{A}$ of
$(\mathcal{C}(X), {\rm sub}_X)$,  the infimum of $\mathcal{A}$ exists; that is  $\sqcap\mathcal{A}=\bigwedge_{C\in \mathcal{C}(X)}\mathcal{A}(C)\to C$.
If $f: X\longrightarrow Y$ is convexity-preserving, then
$$f^\leftarrow:(\mathcal{C}(Y), {\rm sub}_Y)\longrightarrow (\mathcal{C}(X), {\rm sub}_X)$$
 is   infimum-preserving.

\begin{lm}{\rm (\cite[Proposition 2.8]{ZX-Zhang})}\label{lm-equ-cp}
Let $f:X\longrightarrow Y$ be a mapping between two $L$-convex spaces. Then $f$ is convexity-preserving iff
$f^{\rightarrow}(co_X(A))\subseteq co_Y(f^{\rightarrow}(A))$ iff  $co_Y(f^{\rightarrow}(co_X(A)))= co_Y(f^{\rightarrow}(A))$ $(\forall A\in L^{X})$.

\end{lm}


 \begin{lm}{\em (\cite[Proposition  2.10]{Liu-Yue-2024})}\label{liu-yue-pro2.10} Let $X$ be an $L$-convex space. Then $\forall A\in L^X$ and $B\in \mathcal{C}(X)$, we have
${\rm sub}(A, B)={\rm sub}(co(A), B)$.

\end{lm}

\section{Sober $L$-convex spaces}
In \cite{Liu-Yue-2024}, Liu and Yue  introduced a  notion of sober $L$-convex spaces, extending  the study of sober convex spaces to the fuzzy setting. Making use of the inclusion $L$-order between convex sets in
$L$-convex spaces, they  defined algebraic irreducible convex sets in  $L$-convex  spaces and subsequently introduced a notion of sober $L$-convex spaces. In this section, we first recall the concept of
sober $L$-convex spaces defined by Liu and Yue and explore additional properties concerning  this type of     sobriety of $L$-convex spaces.

\begin{dn}{ (\cite[Definition 2.12]{Liu-Yue-2024})}\label{dn-liu-2.12} Let $X$ be an $L$-convex space.
A convex set  $F\in \mathcal{C}(X)$ is  said to be {\em algebraic irreducible} if

{ (1)} $\bigvee_{x\in X}F(x)=1$;

{ (2)} ${\rm sub}(F, \bigvee\nolimits^{\uparrow}_{i\in I}C_{i})=\bigvee\nolimits^{\uparrow}_{i\in I}{\rm sub}(F, C_{i})$
for every directed family $\{C_i\mid i\in I\}\subseteq\mathcal{C}(X)$.
\end{dn}
Let irr$(X)$  denote the collection of all algebraic irreducible convex sets in $X$.
An algebraic irreducible convex set is called a {\it compact convex set} in \cite{Wu-Yao}.

\begin{dn}\label{Liu-Yue-Sober}{ (\cite[Definition 2.13]{Liu-Yue-2024})}
An $L$-convex space   $X$ is said to be {\em sober} if  for each algebraic irreducible convex set  $F\in {\rm irr}(X)$, $F$ is the hull of $1_a$ for a unique $a\in X$,  i.e.,  $F=co(1_a)$.
\end{dn}
 Since we introduced another type of sobriety for $L$-convex spaces in \cite{Wu-Yao}, we provide the following standing assumption to prevent confusion.\\

{\color{black} \noindent {\bf Standing Assumption}. In this paper, the terminology ``sober''  mentioned in this paper always refers to the definition above which was given by Liu and Yue, unless otherwise specified.}\\

Recall that an $L$-convex space $X$ is said to be $S_0$ if for every $x, y\in X$, $C(x)=C(y)$ $(\forall C\in \mathcal{C})$ implies  $x=y$. Using hull operator, Liu and   Yue (\cite[Proposition  2.8]{Liu-Yue-2024}) showed that an $L$-convex space $X$ is  $S_0$ if and only if for all  $x, y\in X$, $co(1_x)=co(1_y)$ implies  $x=y$. Clearly, for every $x\in X$, $co(1_x)$ is an algebraic irreducible convex set. Thus, every sober $L$-convex space is $S_0$.\\

 {\color{black} The following example is the same as  Example 3.4 in \cite{Wu-Yao}. We now   provied a detailed proof that the $L$-convex structure defined below is indeed  a sober $L$-convex spaces in the sense of Definition \ref{dn-liu-2.12}.

\begin{ex}\label{ex-sob-1} Let $L=([0, 1], \otimes)$ be  a commutative  integral quantale with $\otimes$ being the meet operation $\wedge$.  That is to say, $L$ is a frame.
 Define a stratified $L$-convex structure
$$\mathcal{C}=\{a\wedge \phi\mid   a\in [0, 1], \phi: [0, 1]\longrightarrow [0, 1] \text{ is increasing}, \phi\geq id                 \}$$
 on $[0, 1]$, where $id$ sends each $x\in [0,1]$ to $x\in [0,1]$.
 Specifically, a function $\mu: [0, 1]\to [0, 1]$ is a member of  $ \mathcal{C}$ if and only if $\mu$ is an increasing function; and there exists some $a\in [0,1]$ such that
$\mu(x)\geq x$ if $x\in [0, a)$ and $\mu(x)=a$ if $x\in [a, 1]$.

We next show that   in $([0,1],\mathcal{ C})$, every algebraic irreducible convex set is  the hull of $1_x$ for a unique $x\in [0, 1]$.
Let $F\in  {\rm irr}([0,1])$. Since  $\bigvee_{x\in [0,1]}F(x)=1$, we have $F(1)=1$. Write $b=\inf\{x\in [0, 1]\mid F(x)=1\}$ and $A=\{x\in [0, 1]\mid F(x)<1\}$. Clearly, $A=[0, b)$ or $A=[0, b]$. We claim that for every  $x\in A$,  $F(x)=x$. Suppose, for the sake of contradiction,  that there exists $x_0\in A$ such that $F(x_0)\in (x_0, 1)$. For every $n\in \mathbb{N}$, define $F_n:[0, 1]\longrightarrow [0, 1]$ by
$$
\ \ F_n(x)=\left\{\begin{array}{ll}\frac{(n-1)F(x)+x}{n},&x\in A;\\
\ \ \ \ 1,& otherwise.
\end{array}\right.
$$
It is clear that $F_n\in\mathcal{C}$ and $\bigvee_{n\in \mathbb{N}}F_n=F$. Hence, ${\rm sub}(F, \bigvee_{n\in \mathbb{N}}F_n)=1$. For every $n$, we have ${\rm sub}(F, F_n)\leq F(x_0)\rightarrow F_n(x_0)=F_n(x_0)$. Thus,
$$\bigvee_{n\in \mathbb{N}}{\rm sub}(F, F_n)\leq \bigvee_{n\in \mathbb{N}}F_n(x_0)=F(x_0)<1;$$
this is a contradiction. Thus,  for every  $x\in A$,  $F(x)=x$.

We now claim that $b\notin A$, i.e., $F(b)=1$. Assume, for the sake of contradiction, that $b\in A$, i.e., $F(b)=b$.
 For every $n\in \mathbb{N}$, define $G_n:[0, 1]\longrightarrow [0, 1]$ by
$$
\ \ G_n(x)=\left\{\begin{array}{ll} \ \ \ \ x,&x\in [0, b];\\
\frac{n}{2}(x-b)+\frac{b+1}{2},& (b, \frac{1+(n-1)b}{n}];\\
\ \ \ \ 1, & otherwise.

\end{array}\right.$$
It is clear that $G_n\in\mathcal{C}$ and $\bigvee_{n\in \mathbb{N}}G_n=F$. For every $n$, ${\rm sub}(F, G_n)=\frac{b+1}{2}$.
Therefore,
$$\bigvee_{n\in \mathbb{N}}{\rm sub}(F, G_n)=\frac{b+1}{2}.$$ This contradicts to the fact that ${\rm sub}(F, \bigvee_{n\in \mathbb{N}}G_n)=1$. Thus, $b\notin A$. This shows that $F(x)=x$ for all $x\in [0, b)$ and $F(x)=1$ otherwise. It is clear that
$F=co(1_b)$. The uniqueness of $b$ is clear. Thus, $([0,1],\mathcal{ C})$ is a sober $L$-convex space.
\hfill$\Box$
\end{ex}
}
Given an $L$-convex space $X$, for every  $C\in \mathcal{C}(X)$, define an $L$-subset  $\phi(C)$ of $\text{irr}(X)$
by
$\phi(C)(F)={\rm sub}(F, C)\ (\forall F\in \text{irr}(X))$.
Liu and Yue (\cite[Proposition 4.1]{Liu-Yue-2024}) showed that
$\{\phi(C)\mid C\in \mathcal{C}(X)\}$
 forms  an $L$-convex structure on \text{irr}$(X)$. In this paper,  we denote the resulting convex structure
by $\mathcal{C}(\text{irr}(X))$ and  the resulting $L$-convex space $(\text{irr}(X), \mathcal{C}(\text{irr}(X))$) by $\mathbb{S}(X)$. {\color{black}\cite[Proposition 4.2]{Liu-Yue-2024}} states that the $L$-convex space $\mathbb{S}(X)$ is  sober. Liu and Yue  referred to  $\mathbb{S}(X)$ as the {\em sobrification} of $X$.
Here we provide the standard definition of sobrifications corresponding to $L$-convex spaces.
\begin{dn}\label{dn-sobri}
Let $X$ be an  $L$-convex space, let $Z$  be a sober  $L$-convex space and  let $\eta:X\longrightarrow Z$ be a  convexity-preserving mapping.
Then  $Z$ with the mapping $\eta$, or  $Z$ for short,  is called a {\color{black} {\it sobrification}}  of $X$ if for every sober $L$-convex space $Y$ and every  convexity-preserving mapping $f:X\longrightarrow Y$, there exists a unique    convexity-preserving mapping $\overline{f}: Z\longrightarrow Y$ such that $f=\overline{f}\circ \eta$ (see the following figure).
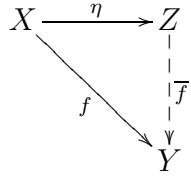
\begin{displaymath}
\xymatrix@=8ex{X\ar[r]^{\eta}\ar[dr]_{f}&Z\ar@{-->}[d]^{\overline{f}}\\  &Y}
\end{displaymath}

\end{dn}
By the  universal property of sobrification, it is easy to see that the sobrification of an  $L$-convex space  is unique up to homeomorphism.

 Let $X$ be an $L$-convex space. Define $\eta_{X}:X\longrightarrow \mathbb{S}(X)$ by $\eta_{X}(x)=co(1_x)$.  By \cite[Lemma 4.3]{Liu-Yue-2024}, $X$ is sober if and only if the corresponding mapping $\eta_{X}$  is a homeomorphism. {\color{black} In this paper, we use $L$-{\bf CS} (resp., $L$-{\bf CS}$_0$, $L$-{\bf SobCS}) to denote the category of $L$-convex spaces (resp., $S_0$ $L$-convex spaces, sober $L$-convex spaces) with convexity-preserving mappings as morphisms.}


Liu and Yue have shown that   $L$-{\bf SobCS} is reflective in $L$-{\bf CS} {\color{black}(see Theorem 5.1 in \cite{Liu-Yue-2024})}. It follows that $\mathbb{S}(X)$, equipped with the mapping $\eta_X$, is a sobrificaion of $X$ in the sense of Definition \ref{dn-sobri}. Therefore,  it is reasonable for Liu and Yue to directly call  $\mathbb{S}(X)$   the sobrificaion of $X$.\\

For an $L$-convex space $X$  and $Y\subseteq X$, define $\mathcal{C}(X)|_{Y}:=\{C|_{Y}\mid C\in \mathcal{C}(X)\}$. The family $\mathcal{C}(X)|_{Y}$ forms  an $L$-convex structure on the background set $Y$. The convex space $(Y, \mathcal{C}(X)|_{Y})$  is referred to as  a {\em subspace} of $(X,  \mathcal{C}(X))$.

For $f:X\longrightarrow Y$, denote $f^{\circ}:X\longrightarrow f(X)$ as the mapping  defined by  $f^{\circ}(x)=f(x)$ for all $x\in X$.
The mapping $f$ is called a {\em subspace embedding} if $f^{\circ}$ is a homeomorphism from $X$ to $f(X)$, where $f(X)$ is  considered as a subspace of $Y$.\\

We present some properties of sober $L$-convex spaces.
\begin{pn}\label{pn-eq-sob} Let $f$ and $g$ be two convexity-preserving mappings from a sober $L$-convex space $X$ to an $S_0$ $L$-convex space $Y$. When $Z=\{x\in X\mid  f(x)=g(x)\}$ is  considered as a subspace of $X$, it is a sober $L$-convex space.
\end{pn}
\noindent{\bf Proof.}
Let $F$ be an algebraic irreducible convex set of $Z$. Define $F^{\ast}\in L^X$ by $F^{\ast}(x)=F(x)$ for $x\in Z$ and $F^{\ast}(x)=0$ for $x\in X-Z$. We claim  that $co_X(F^{\ast})\in {\rm irr}(X)$. Indeed,  for every directed family $\{S_i\mid i\in I\}\subseteq\mathcal{C}(X)$, by $F\in {\rm irr}(Z)$ and Lemma \ref{liu-yue-pro2.10},
\begin{align*}
{\rm sub}_X(co_X(F^{\ast}), \bigvee_{i\in I}\nolimits^{\uparrow}S_i)&={\rm sub}_X(F^{\ast}, \bigvee_{i\in I}\nolimits^{\uparrow}S_i)\\
&={\rm sub}_Z(F, \bigvee_{i\in I}\nolimits^{\uparrow}S_{i}|_{Z})\\
&=\bigvee_{i\in I}\nolimits^{\uparrow}{\rm sub}_Z(F,S_{i}|_{Z})\\
&=\bigvee_{i\in I}\nolimits^{\uparrow}{\rm sub}_X(F^{\ast},S_{i})\\
&=\bigvee_{i\in I}\nolimits^{\uparrow}{\rm sub}_X(co_X(F^{\ast}),S_{i}).
\end{align*}
Since $X$ is sober,  {\color{black} we can let} $x_0\in X$ be the unique element satisfying  $co_X(F^{\ast})=co_X(1_{x_0, X})$. We claim that  $x_0\in Z$; that is,  $f(x_0)=g(x_0)$. Since $Y$ is $S_0$, we only need to show that $co_Y(1_{f(x_0)})=co_Y(1_{g(x_0)})$. For every $E\in \mathcal{C}(Y)$, by $f^{\rightarrow}\dashv f^{\leftarrow}$,
\begin{align*}
1_{f(x_0)}\leq E&\Leftrightarrow f^{\rightarrow}(1_{x_0, X})
\leq E\\
&\Leftrightarrow 1_{x_0, X}\leq f^{\leftarrow}(E)\\
& \Leftrightarrow co_X(F^{\ast})=co_X(1_{x_0, X})\leq f^{\leftarrow}(E)\\
& \Leftrightarrow F^{\ast}\leq f^{\leftarrow}(E)\\
&\Leftrightarrow   f^{\rightarrow} (F^{\ast})\leq E.
\end{align*}
Similarly, $  1_{g(x_0)} \leq E $ if and only if $ g^{\rightarrow} (F^{\ast})\leq E$ for every $E\in \mathcal{C}(Y)$.  Notice that due to the definition of $Z$ and the fact that  $F^{\ast}$  only takes non-zero values on $Z$,  we have $f^{\rightarrow} (F^{\ast}) =g^{\rightarrow} (F^{\ast})$. Altogether, we have $  1_{f(x_0)} \leq E$ if and only if  $  1_{g(x_0)} \leq E$, which shows that $co_Y(1_{f(x_0)})=co_Y(1_{g(x_0)})$. Thus $f(x_0)=g(x_0)$, i.e., $x_0\in Z$. We claim that $co_Z(F)=co_Z(1_{x_0, Z})$.
In fact, for all $C\in \mathcal{C}(X)$,
\begin{align*}
F\leq C|_{Z}&\Leftrightarrow F^{\ast}\leq C\\
&\Leftrightarrow co_X(F^{\ast})\leq C\\
&\Leftrightarrow co_X(1_{x_0, X})\leq C\\
&\Leftrightarrow 1_{x_0, X}\leq C\\
&\Leftrightarrow 1_{x_0, Z}\leq  C|_{Z}.
\end{align*}
This shows that $co_Z(F)=co_Z(1_{x_0, Z})$.
Thus $Z$ is a sober $L$-convex space. \hfill$\Box$\\

\noindent
The preceding discussion  shows that the equalizer of two convexity-preserving mappings from a sober $L$-convex space to an $S_0$ $L$-convex spaces is sober in categorical terms. For the notion of equalizer, please refer to \cite{Category}. Proposition \ref{pn-eq-sob} is
a counterpart of $L$-topological setting (see \cite[Proposition 5.9]{YaoFrm}).\\

Let $X$ and $Y$ be two $L$-convex spaces. The space $Y$ is called a {\em retraction kernel} of $X$ in $L$-{\bf CS} if there exist two convexity-preserving mappings   $r: X\longrightarrow Y$ and $d:Y\longrightarrow X$ such that $r\circ d=id_{Y}$. It is clear that $d$ is an injection and $r$ is a surjection.

%

\begin{pn}\label{pn-ret-sob}  Each retraction kernel of a sober  $L$-convex space  is a sober $L$-convex space.
\end{pn}
\noindent{\bf Proof.}
Let $Y$ be a retraction kernel of $X$ in $L$-{\bf CS}. Then there exist two  convexity-preserving mappings $r: X\longrightarrow Y$ and $d: Y\longrightarrow X$ such that $r\circ d=id_{Y}$. We prove that $Y$ is a sober $L$-convex space. Let $G\in {\rm irr}(Y)$. Then by  \cite[Lemma 4.4]{Liu-Yue-2024},  $co_X( d^{\rightarrow}(G))\in {\rm irr}(X)$. By the sobriety of $X$, there exists a unique $x\in X$ such that $co_X( d^{\rightarrow}(G))=co_X(1_x)$. We claim that $G=co_Y(1_{r(x)})$. On  one hand, by Lemma \ref{lm-equ-cp},
\begin{align*}
 G&=(r\circ d)^{\rightarrow}(G)\\
 &\leq r^{\rightarrow} ( co_{X}(d^{\rightarrow}(G)))\\
 &=  r^{\rightarrow}(co_X(1_x))\\
&\leq   co_Y(1_{r(x)}).
\end{align*}
On the other hand, since $co_X( d^{\rightarrow}(G))=co_X(1_x)$, we have
{\color{black}$$r^{\rightarrow}(co_X( d^{\rightarrow}(G)))(r(x))\geq co_X( d^{\rightarrow}(G))(x)=1.$$}
Thus, $r^{\rightarrow}(co_X( d^{\rightarrow}(G)))(r(x))=1$. Since
{\color{black}\begin{align*}
r^{\rightarrow}(co_X( d^{\rightarrow}(G)))(r(x)
&\leq   co_Y ((r\circ d)^{\rightarrow}(G))(r(x))\\
 &=co_Y(G)(r(x))\\
 &=G(r(x)),
\end{align*}}
we have $G(r(x))=1$, which shows that $co_Y(1_{r(x)})\leq G$.
Thus $G=co_Y(1_{r(x)})$. To verify the uniqueness of $r(x)$, we  need to show that $Y$ is $S_0$. Let $y_1$, $y_2\in Y$ and $y_1\neq y_2$. Since $d$ is an injection, we have $d(y_1)\neq d(y_2)$. Since $X$ is $S_0$, there exists $C\in \mathcal{C}(X)$ such that $ C(d(y_1))\neq C(d(y_2))$. Therefore, $ d^{\leftarrow}(C)(y_1)\neq d^{\leftarrow}(C)(y_2)$. Since $d$ is convexity-preserving, we have $d^{\leftarrow}(C)\in \mathcal{C}(Y)$. This shows that $Y$ is $S_0$.\hfill$\Box$\\


\section{The characterization of sobriety}

Xia in \cite{Xia} {\color{black} introduced} the notions of  quasihomeomorphisms and strict embeddings  in  the framework of convex spaces.  He showed that sober  convex spaces are  precisely injective $S_0$ convex spaces related to strict embeddings.
Following Xia's step, we aim to extend the  notions of  quasihomeomorphisms and strict embeddings to those in the category of $L$-convex spaces. Making use of these notions, we shall provide the characterizations for sobrification and sobriety, respectively.

\begin{dn}\label{dn-str-em} Let  $f:X\longrightarrow Y$ be a  mapping between two   $L$-convex spaces.
 \begin{itemize}
 \item The map $f$ is called a  {\em quasihomeomorphism} if
$f^{\leftarrow}:\mathcal{C}(Y)\longrightarrow  \mathcal{C}(X)$ is a bijection.
\item The map  $f$  is called
a {\em strict embedding} if it is both a quasihomeomorphism and a subspace embedding.
\end{itemize}
\end{dn}
\begin{ex}\label{ex-qusi-inj} Given an  $L$-convex space $X$, $\eta_{X}:X\longrightarrow \mathbb{S}(X)$ is a quasihomeomorphism.
It is routine to check that
for every $C\in \mathcal{C}(X)$,  $\eta_{X}^{\leftarrow}(\phi(C))=C$ and $\phi\circ \eta_{X}^{\leftarrow}(\phi(C))=\phi(C)$. Thus, $\eta_{X}^{\leftarrow}\circ \phi=id_{\mathcal{C}(X)}$ and $\phi\circ\eta_{X}^{\leftarrow}=id_{\mathcal{C}({\rm irr}(X))}$.
\end{ex}

The following shows that  in the category of $S_0$ $L$-convex spaces,  quasihomeomorphisms  coincide with strict embeddings.

 \begin{pn}Let $f:X\longrightarrow Y$ be a  mapping between two $S_0$ $L$-convex spaces. Then  $f$ is a subspace embedding if and only if $f^{\leftarrow}:\mathcal{C}(Y)\longrightarrow  \mathcal{C}(X)$ is a surjection.
\end{pn}
\noindent{\bf Proof.} {\it Necessity.} For every $C\in \mathcal{C}(X)$, since $f$ is a subspace embedding, we have $(f^{\circ})^{\rightarrow}(C)\in \mathcal{C}(Y)|_{f(X)}(:=\{D|_{f(X)}\mid D\in \mathcal{C}(Y)\})$. So, there exits $D\in \mathcal{C}(Y)$ such that $(f^{\circ})^{\rightarrow}(C)=D|_{f(X)}$.  We  prove that $f^{\leftarrow}(D)=C$. For every $x\in X$, since $f$ is an injection, we have
\begin{align*}
f^{\leftarrow}(D)(x)&={\color{black}D(f(x))}\\
&=D|_{f(X)}(f(x))\\
&=(f^{\circ})^{\rightarrow}(C)(f(x))\\
&=C(x).
\end{align*}
Thus $f^{\leftarrow}(D)=C$. By the arbitrariness of  $C\in \mathcal{C}(X)$, we have that $f^{\leftarrow}:\mathcal{C}(Y)\longrightarrow  \mathcal{C}(X)$ is a surjection.

{\it Sufficiency.} First, we show that $f$ is an injection. Let $x, y\in X$ and $x\neq y$. Since $X$ is $S_0$, there exists $C\in \mathcal{C}(X)$ such that $C(x)\neq C(y)$. As $f^{\leftarrow}$ is a surjection, there exists $E\in \mathcal{C}(Y)$ such that $f^{\leftarrow}(E)=C$. Therefore
$$E(f(x))=f^{\leftarrow}(E)(x)\neq f^{\leftarrow}(E)(y)=E(f(y)).$$ Thus $f(x)\neq f(y)$, showing that $f$ is an injection.

Next, we show that $f^{\circ}:(X, \mathcal{C}(X))\longrightarrow (f(X), \mathcal{C}(Y)|_{f(X)})$ is a  convex-to-convex mapping. Let $C\in \mathcal{C}(X)$. Since $f^{\leftarrow}$ is a surjection, there exits $E\in \mathcal{C}(Y)$ such that $f^{\leftarrow}(E)=C$. We  claim that $(f^{\circ})^{\rightarrow}(C)=E|_{f(X)}$. In fact, for every $y\in f(X)$, assume $y=f(x)$ for some $x\in X$,
\begin{align*}
(f^{\circ})^{\rightarrow}(C)(y)&=(f^{\circ})^{\rightarrow}(f^{\leftarrow}(E))(y)\\
&=(f^{\circ})^{\rightarrow}(f^{\leftarrow}(E))(f(x))=f^{\leftarrow}(E)(x)\\
&=E(f(x))=E(y).
\end{align*}
 Thus, $(f^{\circ})^{\rightarrow}(C)=E|_{f(X)}$, showing that $f^{\circ}:(X, \mathcal{C}(X))\longrightarrow \mathcal{C}(Y)|_{f(X)})$ is  convex-to-convex. In conclusion, $f$ is a subspace embedding.\hfill$\Box$\\

Now we obtain the characterization for sobrification. The $L$-ordered {\color{black} structure} on the family of convex sets plays a crucial role in this characterization. We first see a lemma.
\begin{lm}\label{lm-bij-iso} Let $f:P\longrightarrow Q$ be a bijection between two complete $L$-ordered sets. Then  $f$ is  supremum-preserving (or, infimum-preserving) if and only if $f$ is an  $L$-order-isomorphism.
\end{lm}
\noindent {\bf Proof.} We  assume that $f$ is  supremum-preserving.  The proof of the case that $f$ is  infimum-preserving is similar.

{\it {\color{black}Necessity}}.
For every $x, y\in P$, we have $f(y)=f(\sqcup{\downarrow} y)=\sqcup f^{\rightarrow}({\downarrow} y)$. By Condition (1) in Definition \ref{dn-sup-inf},
$$e_P(x, y)=f^{\rightarrow}({\downarrow} y)(f(x))\leq e_Q(f(x), f(y)).$$
On the other hand,
$$f(\sqcup f^{\leftarrow} ({\downarrow} f(y)))= \sqcup  f^{\rightarrow}( f^{\leftarrow} ({\downarrow} f(y)))=\sqcup{\downarrow} f(y)=f(y). $$
Since $f$ is a bijection, $\sqcup f^{\leftarrow} ({\downarrow} f(y))=y$. By Condition (1) in Definition \ref{dn-sup-inf},
$$f^{\leftarrow} ({\downarrow} f(y))(x)=e_Q(f(x), f(y))\leq e_P(x, y).$$ Thus
$e_P(x, y)=e_Q(f(x),f(y))$, which shows that $f$ is  an $L$-order-isomorphism.

{\it {\color{black}Sufficiency}}. Let $g$ be the inverse mapping of $f$. For every $A\in L^P$ and $x\in Q$,
 \begin{align*}
{\rm sub}(f^{\rightarrow}(A), {\downarrow} x)&=\bigwedge_{y\in Q} f^{\rightarrow}(A)(y) \to e(y,x)\\
&=\bigwedge_{t\in P}A(t)\to e(f(t),x)\\
&=\bigwedge_{t\in P}A(t)\to e(t,g(x))\\
&={\rm sub}(A, {\downarrow} g(x))\\
&=e_P(\sqcup A, g(x))\\
&=e_Q(f(\sqcup A), x).
\end{align*}
Thus $f(\sqcup A)=\sqcup f^{\rightarrow}(A)$. This shows that $f$ is a supremum-preserving mapping.
\hfill$\Box$\\

\begin{tm}\label{tm-char-sobri}{\em (Characterization {\color{black} Theorem I}: for sobrification)}
Let $X$  and $Y$ be two $L$-convex spaces.  If $Y$ is   sober, then the following  are equivalent:

$(1)$ $Y$ is a sobrification of $X$;

$(2)$ there exits a quasihomeomorphism  from $X$ to $Y$;

$(3)$  $(\mathcal{C}(X), {\rm sub}_X)$  and $(\mathcal{C}(Y), {\rm sub}_Y)$ are $L$-order-isomorphic.
\end{tm}
\noindent {\bf Proof.}
$(1)\Rightarrow (2)$: Since  $Y$ and  $\mathbb{S}(X)$ are  sobrifications of $X$, there exists a homeomorphism $f:\mathbb{S}(X)\longrightarrow Y$. By Example \ref{ex-qusi-inj},  $\eta_{X}:X\longrightarrow \mathbb{S}(X)$ is a quasihomeomorphism. Since a composition of quasihomeomorphism and homeomorphism is a quasihomeomorphism, it follows that $f\circ\eta_{X}$ is a quasihomeomorphism from $X$ to $Y$.

$(2)\Rightarrow (3)$: It is clear that $j^{\leftarrow}: (\mathcal{C}(Y), {\rm sub}_Y)\to(\mathcal{C}(X), {\rm sub}_X)$ is an infimum-preserving  bijection. By Lemma \ref{lm-bij-iso}, we have $j^{\leftarrow}$ is an $L$-order-isomorphism.

$(3)\Rightarrow(1)$: For an $L$-convex space $X$, it is easy to see that  ${\rm irr}(X)$  and  $\mathcal{C}({\rm irr}(X))$  are defined by the corresponding $L$-ordered structure on $(\mathcal{C}(X), {\rm sub}_X)$. Therefore,
$$({\rm irr}(Y), \mathcal{C}({\rm irr}(Y)))\cong({\rm irr}(X), \mathcal{C}({\rm irr}(X))).$$
Since $Y$ is  sober,  it follows that $(Y, \mathcal{C}(Y))\cong({\rm irr}(Y), \mathcal{C}({\rm irr}(Y)))$. Thus
$$(Y, \mathcal{C}(Y))\cong({\rm irr}(X), \mathcal{C}({\rm irr}(X))).$$
Note that   ${\rm irr}(X)$ is a sobrification of $X$. Hence,
 $Y$ is also a sobrification of $X$.

\hfill$\Box$\\

{\color{black}The following example  shows
that Theorem \ref{tm-char-sobri}(3) can be effectively used to identify the sobrification of
an $L$-convex space.

\begin{ex}\label{ex-sobri} Let $L=([0, 1], \otimes)$ be   a commutative  integral quantale with $\otimes$ being $\wedge$.     Define a convex structure on $[0,1)$ as follows:
$$\mathcal{C}^{\prime}=\{a\wedge \psi\mid   a\in [0, 1], \psi: [0, 1)\longrightarrow [0, 1] \text{ is increasing}, \psi\geq id                 \},$$
where $id$ sends each $x\in [0,1)$ to $x\in [0,1]$. Clearly, $([0,1), \mathcal{C}^{\prime})$ is a subspace of   $([0,1], \mathcal{C})$  defined in Example  \ref{ex-sob-1}.  We now show that  $id: [0, 1)\longrightarrow [0,1]$ is an algebraic irreducible convex set. Obviously,  we have $id\in \mathcal{C}$ and $\bigvee_{x\in [0,1)}id(x)=1$. For every $\mu\in \mathcal{C}$,  there exists some $a\in [0,1]$ such that
$\mu(x)\geq x$ if $x\in [0, a)$ and $\mu(x)=a$ if $x\in [a, 1)$. Thus, we have
 $${\rm sub}(id, \mu)=\bigvee_{x\in[0, 1)}\mu(x)=a.$$ For a directed subset $\{\mu_i\mid i\in I\}\subseteq\mathcal{C}$, we have
$${\rm sub}(id, \bigvee\nolimits^{\uparrow}_{i\in I}\mu_i)=\bigvee_{x\in[0, 1)}\bigvee\nolimits^{\uparrow}_{i\in I}\mu_i(x)=\bigvee\nolimits^{\uparrow}_{i\in I}\bigvee_{x\in[0, 1)}\mu_i(x)=\bigvee\nolimits^{\uparrow}_{i\in I}{\rm sub}(id, \mu_i).
$$
 Thus,  $id$ is an algebraic irreducible convex set. However, $id$ is not the hull of  $1_x$ for a unique $x\in [0, 1)$. Therefore, $([0, 1), \mathcal{C}^{\prime})$ is not sober. It is straightforward to check  that $(\mathcal{C},  {\rm sub}_{[0,1]})\cong(\mathcal{C}^{\prime}, {\rm sub}_{[0,1)})$. By Theorem \ref{tm-char-sobri} (3),  the sober $L$-convex space $([0, 1], \mathcal{C})$ is  a sobrification of $([0, 1), \mathcal{C}^{\prime})$.
\hfill$\Box$
\end{ex}

\begin{rk}\label{rk-sob}
In Example \ref{ex-sobri}, similar to Example 3.4 in \cite{Wu-Yao}, we can deduce  that $A$ is nonempty finite $L$-subset of $[0, 1)$ if and only if   there exists a nonempty finite subset $K\subseteq_{fin} [0,1)$ such that $A =\chi_{K}$. For the definition of finite $L$-subset, please refer to \cite[Definition 3.1]{Wu-Yao}.
For $K\subseteq_{fin} [0,1)$, write ${\rm min} (K)=b_0$. It is clear that  $co(\chi_{K})=co(1_{b_0})$. Thus, $\mathcal{C}$ is sober in the sense of Wu and Yao (see  Definition 3.1 in \cite{Wu-Yao}). However, Example \ref{ex-sobri} shows that $\mathcal{C}$ is not sober in the sense of Liu and Yue (see Definition \ref{Liu-Yue-Sober}). Thus, we here have provided an answer to the unresolved problem propopsed in  \cite{Wu-Yao}. To prevent confusion, we can rename the sobriety defined in \cite{Wu-Yao} as {\it  ``weak sobriety''}.
\end{rk}}

For every  convexity-preserving mapping $f:X\longrightarrow Y$, by  \cite[Lemma 4.4]{Liu-Yue-2024}, we can define  a mapping
$$\mathbb{S}(f):\mathbb{S}(X)\longrightarrow \mathbb{S}(Y)$$
 by $\mathbb{S}(f)(F)=co_Y(f^{\rightarrow}(F))$ for every $F\in {\rm irr}(X)$.

\begin{lm}\label{lm-cp-cp} {\em (1)} The mapping $\mathbb{S}(f):\mathbb{S}(X)\longrightarrow \mathbb{S}(Y)$ is convexity-preserving;

 {\em (2)} The equality $\mathbb{S}(f)\circ\eta_{X}=\eta_{Y}\circ f$ holds.
\end{lm}
\noindent{\bf Proof.}
 (1) For every $E$ in  $\mathcal{C}(Y)$ and  $G \in \text{irr}(Y)$, it holds that
{\color{black} \begin{align*}
\mathbb{S}(f)^{\leftarrow}(\phi(E))(G)
&=\phi(E)( co_{Y}(f^{\rightarrow}(G))\\
&={\rm sub}_Y(co_{Y}(f^{\rightarrow}(G)), E)\\
&={\rm sub}_Y(f^{\rightarrow}(G), E)\\
&={\rm sub}_X(G, f^{\leftarrow}(E))\\
&=\phi(f^{\leftarrow}(E))(G).
\end{align*}}
Thus, $\mathbb{S}(f)^{\leftarrow}(\phi(E))=\phi(f^{\leftarrow}(E))\in \mathcal{C}({\rm irr}(X))$, which shows that $\mathbb{S}(f)$ is convexity-preserving.

(2) The equality $\mathbb{S}(f)\circ\eta_{X}=\eta_{Y}\circ f$ holds by Lemma \ref{lm-equ-cp}.
\hfill$\Box$\\

Let $f:X\longrightarrow Y$ be a mapping, with  $X^{\prime}\subseteq X$ and $Y^{\prime}\subseteq Y$. If $f(X^{\prime})\subseteq Y^{\prime}$,
 we  denote $f|_{X^{\prime}}^{Y^{\prime}}:X^{\prime}\longrightarrow Y^{\prime}$ as the mapping  defined by $f|_{X^{\prime}}^{Y^{\prime}}(x)=f(x)$ for all $x\in X^{\prime}$.

\begin{pn}\label{pn-flat-homo} Let $X$ and  $Y$ be two $L$-convex spaces, and let $f:X\longrightarrow Y$ be a convexity-preserving mapping. Then following  are equivalent:

$(1)$  $f$ is a quasihomeomorphism;

$(2)$   $\mathbb{S}(f): \mathbb{S}(X)\longrightarrow \mathbb{S}(Y)$  is a homeomorphism.

\end{pn}
\noindent {\bf Proof.} $(1)\Rightarrow (2)$: Let $f_{\ast}: \mathcal{C}(X)\longrightarrow \mathcal{C}(Y)$ be the inverse mapping of $f^{\leftarrow}: \mathcal{C}(Y)\longrightarrow\mathcal{C}(X)$. Clearly, $f_{\ast}$ is an $L$-order-isomorphism. Notice that the notion of algebraic irreducible convex sets  is given by means of the inclusion  $L$-order between  convex sets. So,  $f_{\ast}|_{{\rm irr}(X)}^{{\rm irr}(Y)}:{\rm irr}(X)\longrightarrow {\rm irr}(Y)$ is an $L$-order-isomorphism, whose inverse mapping is $f^{\leftarrow}|_{{\rm irr}(Y)}^{{\rm irr}(X)}:{\rm irr}(Y)\longrightarrow{\rm irr}(X)$. We claim that $f_{\ast}|_{{\rm irr}(X)}^{{\rm irr}(Y)}=\mathbb{S}(f)$; that is,  for all $F\in {\rm irr}(X)$, $f_{\ast}(F)=\mathbb{S}(f)(F)(=co_Y(f^{\rightarrow}(F)))$. In fact,  for every $D\in  \mathcal{C}(Y)$, we have
$$f^{\rightarrow}(F)\leq D\Longleftrightarrow F\leq f^{\leftarrow}(D)
\Longleftrightarrow  f_{\ast}(F)\leq  D$$
Notice that  $f_{\ast}(F)\in \mathcal{C}(Y)$. We have $co_Y(f^{\rightarrow}(F))=f_{\ast}(F)$.
 By Lemma \ref{lm-cp-cp}, $\mathbb{S}(f)=f_{\ast}|_{{\rm irr}(X)}^{{\rm irr}(Y)}:{\rm irr}(X)\longrightarrow {\rm irr}(Y)$ is convexity-preserving. We now only need to show that $\mathbb{S}(f)$ is convex to convex.
For every $C\in \mathcal{C}(X)$ and $G\in {\rm irr}(Y)$,
{\color{black}\begin{align*}
(\mathbb{S}(f))^{\rightarrow}(\phi(C))(G)
&=\phi(C)((\mathbb{S}(f))^{-1}(G))\\
&=\phi(C)(f^{\leftarrow}(G))\\
&={\rm sub}_X (f^{\leftarrow}(G), C)\\
&={\rm sub}_Y(G, f_\ast(C))\\
&=\phi(f_\ast(C))(G).
\end{align*}}
Thus, $(\mathbb{S}(f))^{\rightarrow}(\phi(C))=\phi(f_\ast(C))$. This shows that $\mathbb{S}(f):{\rm irr}(X)\longrightarrow {\rm irr}(Y)$ is   convex to convex,  hence  a homeomorphism.

$(2)\Rightarrow (1)$: Let $\mathbb{S}(f)_{*}$ denote the inverse mapping of $\mathbb{S}(f)$.
Define  $f_{\ast}:\mathcal{C}(X)\longrightarrow\mathcal{C}(Y)$ as the composition $(\eta_{Y})^{\leftarrow}\circ((\mathbb{S}(f)_{*})^{\leftarrow}\circ\phi$ (see figure below).
\begin{displaymath}
\xymatrix@=12ex{\mathcal{C}(\text{irr}(X))\ar[r]^{(\mathbb{S}(f)_{*})^{\leftarrow}}&\mathcal{C}(\text{irr}(Y))\ar[d]^{(\eta_{Y})^{\leftarrow}}\\ \mathcal{C}(X)\ar[r]^{f_{\ast}}\ar[u]^{\phi}&\mathcal{C}(Y)}
\end{displaymath}
\vskip 3pt
 \noindent For  every $C\in \mathcal{C}(X)$, by Lemma \ref{lm-cp-cp}(2), we have
{\color{black}\begin{align*}
f^{\leftarrow}\circ f_{\ast}(C)
&=f^{\leftarrow}\circ (\eta_{Y})^{\leftarrow}\circ((\mathbb{S}(f)_{*})^{\leftarrow}\circ\phi(C)\\
&=(\eta_{Y}\circ f)^{\leftarrow}\circ ((\mathbb{S}(f)_{*})^{\leftarrow}(\phi(C))\\
&=(\mathbb{S}(f)\circ \eta_{X})^{\leftarrow}\circ ((\mathbb{S}(f)_{*})^{\leftarrow}(\phi(C))\\
&=\eta_{X}^{\leftarrow}\circ \mathbb{S}(f)^{\leftarrow}\circ ((\mathbb{S}(f)_{*})^{\leftarrow}(\phi(C))\\
&=\eta_{X}^{\leftarrow}\circ ((\mathbb{S}(f)_{*}\circ \mathbb{S}(f))^{\leftarrow}(\phi(C))\\
&=\eta_{X}^{\leftarrow}(\phi(C))\\
&=C.
\end{align*}}
Therefore, $f^{\leftarrow}\circ f_{\ast}=id_{\mathcal{C}(X)}$.\\

For  every $D\in \mathcal{C}(Y)$, by Lemma \ref{lm-cp-cp}(1), we have
{\color{black}\begin{align*}
 f_{\ast}\circ f^{\leftarrow}(D)
&=(\eta_{Y})^{\leftarrow}\circ ((\mathbb{S}(f)_{*})^{\leftarrow}\circ\phi(f^{\leftarrow}(D))\\
&=(\eta_{Y})^{\leftarrow}\circ((\mathbb{S}(f)_{*})^{\leftarrow}(\mathbb{S}(f))^{\leftarrow}(\phi(D))&\\
&=(\eta_{Y})^{\leftarrow}(\mathbb{S}(f)\circ \mathbb{S}(f)_{*})^{\leftarrow}(\phi(D))\\
&=(\eta_{Y})^{\leftarrow}(\phi(D))\\
&=D.
\end{align*}}
Therefore, $f_{\ast}\circ f^{\leftarrow}=id_{\mathcal{C}(Y)}$. Thus, $f^{\leftarrow}:\mathcal{C}(Y)\longrightarrow \mathcal{C}(X)$ is a bijection; hence $f$ is a quasihomeomorphism.\hfill$\Box$\\

For a category $\mathcal{C}$, we use the symbol $Mor(\mathcal{C})$ to denote the class of morphisms of $\mathcal{C}$. Readers can refer to
 \cite[Section II-3]{Gierz} for the notion of a {\rm $J$-injective} object in the terminology of  category. For  convenience, we list it below.
 \begin{dn}
{ (\cite{Gierz})} Let $\mathcal{M}\subseteq Mor(\mathcal{C})$ be a class of monomorphisms which is closed under multiplication with isomorphisms, i.e., for every $j\in \mathcal{M}$ and every isomorphism $i$, ones has $j\circ i\in \mathcal{M}$ and $i\circ j \in \mathcal{M}$. Then an object $X$ is said to be {\em $\mathcal{M}$-injective} if for every mapping $j: Y\longrightarrow Z$ in $\mathcal{M}$ and every morphism $f: Y\longrightarrow X$, there exists a morphism $g: Z\longrightarrow X$ with $f=g\circ j$; that is the following diagram commutes.
\begin{displaymath}
\xymatrix@=8ex{Y\ar[r]^{j}\ar[d]_{f}&Z\ar@{-->}[dl]^{g}\\  X}
\end{displaymath}

\end{dn}

In $L$-$\mathbf{CS}_0$, let $\mathcal{M}_1$  and $\mathcal{M}_2$ be the classes of all subspace embeddings and  all strict embeddings, respectively. An $S_0$ $L$-topological space $X$ is said to be {\it  injective} (resp., {\it strictly injective}) if  $X$  is $\mathcal{M}_1$-injective (resp.,   $\mathcal{M}_2$-injective) in the category $L$-$\mathbf{CS}_0$.

We now provide the characterization of  sober $L$-convex space in $L$-{\bf CS}$_0$.

\begin{tm}\label{tm-inj-sob}{\rm (Characterization {\color{black} Theorem II}: for sober $L$-convex space)} The $L$-convex  space  $X$ is a strictly injective $S_0$ $L$-convex space if and only if $X$ is a  sober $L$-convex space.
\end{tm}
\noindent{\bf Proof.}
{\it {\color{black} Necessity}.} It follows from Example \ref{ex-qusi-inj} that  $\eta_{X}$  is a strict embedding.  As $X$ is strictly injective,  we obtain a convexity-preserving mapping $r_X: \mathbb{S}(X)\longrightarrow X$ such that $r_X\circ \eta_X=id_X$.
It follows from Proposition \ref{pn-ret-sob} that $X$ is also a sober $L$-convex space.

{\it {\color{black} Sufficiency}.} Let $X$ be a sober $L$-convex space, and     let $j: Y\longrightarrow Z$ be a strict embedding between two  $S_0$ $L$-convex spaces. Then $\eta_{X}:X\longrightarrow \mathbb{S}(X)$ is a bijection.  By Proposition \ref{pn-flat-homo},  $\mathbb{S}(j):\mathbb{S}(Y)\longrightarrow \mathbb{S}(Z)$ is a homeomorphism.  For every convexity-preserving mapping $f: Y\longrightarrow X$,
 define $\overline{f}:Z\longrightarrow X$ as the composition $(\eta_{X})^{-1}\circ \mathbb{S}(f)\circ (\mathbb{S}(j))^{-1}\circ \eta_{Z}$ (see figure below).
\begin{displaymath}
\xymatrix@=8ex{\mathbb{S}(Z)\ar[r]^{{\color{black} (\mathbb{S}(j))^{-1}}}&\mathbb{S}(Y)\ar[r]^{\mathbb{S}(f)}&\mathbb{S}(X)\ar[d]^{(\eta_{X})^{-1}}\\ Z\ar[rr]^{\overline{f}}\ar[u]^{\eta_{Z}}&&X}
\end{displaymath}
\noindent By Lemma \ref{lm-cp-cp}(2),  we have
{\color{black}\begin{align*}
\overline{f}\circ j
&=(\eta_{X})^{-1}\circ \mathbb{S}(f)\circ (\mathbb{S}(j))^{-1}\circ \eta_{Z}\circ j\\
&=(\eta_{X})^{-1}\circ \mathbb{S}(f)\circ  (\mathbb{S}(j))^{-1}\circ \mathbb{S}(j)\circ \eta_{Y}\\
&=(\eta_{X})^{-1}\circ \mathbb{S}(f)\circ \eta_{Y}\\
&=(\eta_{X})^{-1}\circ \eta_{X}\circ f=f.
\end{align*}}
Thus, $\overline{f}\circ j=f$.
Therefore, $X$ a strictly injective $S_0$ $L$-convex space.\hfill$\Box$

\begin{rk} It is straightforward to verify  that the  mapping $\overline{f}$ constructed in the above proof  is unique.
 In fact,
if $g:Z\longrightarrow X$  satisfies $f=g\circ j$, then by Lemma \ref{lm-cp-cp}(2), we have
\begin{align*}
\overline{f}&=(\eta_{X})^{-1}\circ \mathbb{S}(f)\circ (\mathbb{S}(j))^{-1}\circ \eta_{Z}\\
&=(\eta_{X})^{-1}\circ \mathbb{S}(g\circ j)\circ (\mathbb{S}(j))^{-1}\circ \eta_{Z}\\
&=(\eta_{X})^{-1}\circ \mathbb{S}(g)\circ \mathbb{S}(j)\circ (\mathbb{S}(j))^{-1}\circ \eta_{Z}\\
&=(\eta_{X})^{-1}\circ \mathbb{S}(g)\circ \eta_{Z}\\
&=(\eta_{X})^{-1}\circ \eta_{X}\circ g\\
&=id_X\circ g=g.
\end{align*}
\end{rk}
\section{Conclusion remarks}
{\color{black} The primary focus of this paper is to characterize  sober $L$-convex spaces with a commutative integral quantale $L$ as the truth value table.
By means of the fuzzy ordered methods, we have  derived a characterization theorem  for sobrification (see Theorem \ref{tm-char-sobri}).
 From a categorical perspective, we have  proved that an $L$-convex  space  $X$ is  sober  if and only if it is a strictly injective $S_0$ $L$-convex space. Moreover, we have solved a problem left in \cite{Wu-Yao} by presenting a  counterexample; specifically,   the sobriety defined by Wu and Yao \cite{Wu-Yao} does not imply the sobriety defined  by Liu and Yue
 \cite{Liu-Yue-2024} (see Remark \ref{rk-sob}). Consequently, the fuzzy ordered methods and categorical methods can be effectively combined  in the study  of fuzzy convex theory.

The methods of studying the sobriety of (fuzzy) topological spaces in \cite{Gierz,Noor,Zhao-Zhang} can also be applied in studying   the sobriety of (fuzzy) convex spaces. We have found that many results in \cite{Gierz,Noor,Zhao-Zhang}  also hold in the framework of (fuzzy) convex spaces.  This raises a natural question: why does this phenomenon occur?  In the future, we will investigate the fundamental principles behind this phenomenon, abstract out a unified method, and further explore the categorical properties of sobriety in more general closure spaces as presented  in \cite{Erne-stone} and \cite{ZX-Zhang}.}

\vspace*{.15in}

\noindent{\bf Acknowledgements.}

This paper is supported by National Natural Science Foundation of China
(12231007, 12371462), Jiangsu Provincial Innovative and Entrepreneurial
Talent Support Plan (JSSCRC2021521).


\vspace*{.15in}



\begin{thebibliography}{99}
	
%
\bibitem{Category} J. Ad\'{a}mek, H. Herrlich, G. E. Strecker, {\it Abstract and Concrete Categories}, Wiley, New York, 1990.
\bibitem{Erne-stone} M. Ern\'{e}, {\it General Stone duality}, Topol. Appl., {\bf 137} (2004), 125--158. https://doi.org/10.1016/S0166-8641(03)00204-9
%
\bibitem{M.H.Es}M. H. Escard\'{o}, R. C. Flagg, {\it Semantic Domains, injective spaces and monad: extended abstract},  Electron. Notes Theor. Comput. Sci., {\bf 20} (1999), 229--244. https://doi.org/10.1016/S1571-0661(04)80077-X
\bibitem{Varlet} S. P. Franklin, {\it Some results on order convexity}, Amer. Math. Monthly,  {\bf 69} (1962), 357--359. https://doi.org/10.1080/00029890.1962.11989897

\bibitem{Jean-2013} J. Goubault-Larrecq, {\it Non-Hausdorff Topology and Domain Theory}, Cambridge University Press,  Cambridge, 2013.
\bibitem{Gierz} G. Gierz, K. H. Hofmann, K. Keimel, J. D. Lawson, M. Mislove, D. S. Scott, {\it Continuous Lattices and Domains},  Cambridge University Press,  Cambridge, 2003.



\bibitem{MVTop}U. H\"{o}hle, {\it Many-valued Topology and Its Applications}, Kluwer Academic Publishers, Boston, Dordrecht, London, 2001.


\bibitem{topology2} H. Komiya, {\it Convexity on a topological space}, Fund. Math.,  {\bf 111} (1981), 107--113. https://doi.org/10.4064/fm-111-2-107-113
\bibitem{Liu-Meng-Yue}
M. Y. Liu, Y. L. Yue, X. W. Wei, {\it Frame-valued Scott open set monad and its algebras}, Fuzzy Sets Syst.,
{\bf 460} (2023), 52--71.\\
https://doi.org/10.1016/j.fss.2022.11.002
\bibitem{Liu-Yue-2024} M. Y. Liu, Y. L. Yue,  {\it The reflectivity of the category of stratified $L$-algebraic closure spaces}, Iran. J. Fuzzy Syst.,  {\bf  21} (2024) 117--127.

\bibitem{Algebra}E. Marczewski, {\it Independence in abstract algebras results and problems}, Colloq. Math.,   {\bf 14} (1966), 169--188.  https://doi.org/10.4064/cm-14-1-169-188
\bibitem{Mar-free-2010} Y. Maruyama, {\it Fundamental results for pointfree convex geometry}, Ann. Pure Appl. Logic, {\bf 161} (2010),
1486--1501.\\
 https://doi.org/10.1016/j.apal.2010.05.002
\bibitem{Maruyama} Y. Maruyama, {\it Lattice-valued fuzzy convex geometry}, Comput. Geom. Discrete Math., {\bf 164} (2009), 22--37.


\bibitem{algebra4}J. Nieminen, {\it The ideal structure of simple ternary algebras}, Colloq. Math., {\bf 40} (1978), 23--29. https://doi.org/10.4064/cm-40-1-23-29
\bibitem{Noor} R. Noor, A. K. Srivastava, The categories $L$-{\bf Top}$_0$ and $L$-{\bf Sob} as epireflective hulls, Soft Comput., {\bf 18} (2014),
1865--1871. https://doi.org/10.1007/s00500-014-1315-8
\bibitem{Pang-2023} B. Pang,  {\it Quantale-valued convex structures as lax algebras}, Fuzzy Sets Syst.,  {\bf 473} (2023), 108737. https://doi.org/10.1016/j.fss.2023.108737
\bibitem{PangShi}B. Pang, F. G. Shi, {\it Subcategories of the category of $L$-convex spaces}, Fuzzy Sets Syst., {\bf 313} (2017), 61--74.   https://doi.org/10.1016/j.fss.2016.02.014


\bibitem{Pang2}B. Pang,  F. G. Shi, {\it Fuzzy counterparts of hull operators and interval operators in the framework of $L$-convex spaces},
Fuzzy Sets Syst., {\bf 369} (2019), 20--39. https://doi.org/10.1016/j.fss.2018.05.012






\bibitem{RosaDc} M. V. Rosa, {\it A study of fuzzy convexity with special reference to separation properties}, Cochin University of science and Technology, Cochin, India, 1994.


\bibitem{Quantale} K. I. Rosenthal, {\it Quantales and Their Applications},
Longman House, Burnt Mill, Harlow, 1990.

\bibitem{Su-Li-Intel}S. H. Su, Q. G. Li, {\it The category of algebraic $L$-closure systems},  J. Intell. Fuzzy Syst. {\bf 33} (2017), 2199--2210. https://doi.org/10.3233/JIFS-16452
\bibitem{Sun-Pang-2024}L. C. Sun, B. Pang, {\it Dual equivalence between the categories of sober convex spaces
and algebraic lattices}, 2024, preprint.
\bibitem{Shen1}C. Shen, F. G. Shi, {\it Characterizations of L-convex spaces via domain theory}, Fuzzy Sets Syst., {\bf 380} (2020), 44--63.
https://doi.org/10.1016/j.fss.2019.02.009
\bibitem{Shen}C. Shen, S. J. Yang, D. S. Zhao, F.G. Shi, {\it Lattice-equivalence of convex spaces}, Algebra Univer. {\bf 80} (2019),  26.  https://doi.org/10.1007/s00012-019-0600-x








\bibitem{Van1984} M. Van De Vel, {\it Binary convexities and distributive lattices}, Proc.  London  Math. Soc., {\bf 48} (1984), 1--33.   https://doi.org/10.1112/plms/s3-48.1.1


\bibitem{topology1} M. Van De Vel, {\it On the rank of a topological convexity}, Fund. Math.,  {\bf 119} (1984), 17--48. https://doi.org/10.4064/fm-119-2-101-132

\bibitem{Vanbook} M. Van De Vel, {\it Theory of convex spaces}, North-Holland, Amsterdam, 1993.
\bibitem{Wu-Guo-Li} M. Y. Wu, L. K. Guo, Q. G. Li, {\it New representatons of algebraic domains and algebraic L-domains via closure systems}, Semigroup Forum,
{\bf 103} (2021), 700--712. https://doi.org/10.1007/s00233-021-10209-7
%

\bibitem{Wu-Yao} G. J. Wu, W. Yao, {\it  Sober $L$-convex spaces and $L$-join-semilattices},  Iran. J. Fuzzy Syst., online.
\bibitem{Xia}C. C. Xia, {\it Some further results on pointfree convex geometry}, Algebra Univers., {\bf 85} (2024)
 20. https://doi.org/10.1007/s00012-024-00847-7
\bibitem{XU-mao-form} L. S. Xu,  X. X. Mao, {\it Formal topological characterizations of various continuous domains}, Comput. Math. Appl. {\it 56} (2008) 444--452.\\ https://doi.org/10.1016/j.camwa.2007.12.010
\bibitem{PartI}W. Yao, {\it Quantitative domains via fuzzy sets: Part I: Continuity
of fuzzy directed-complete poset}, Fuzzy Sets Syst., {\bf 161} (2010), 983--987. https://doi.org/10.1016/j.fss.2009.06.018
\bibitem{YaoFrm}W. Yao, {\it A survey of fuzzifications of frames, the Papert-Papert-Isbell adjunction and sobriety},
Fuzzy Sets Syst., {\bf 190} (2012), 63--81.\\
 https://doi.org/10.1016/j.fss.2011.07.013
\bibitem{Yao2016}W. Yao, S.-E. Han, {\it A Stone-type duality for sT$_{0}$ stratified Alexandrov $L$-topological
spaces}, Fuzzy Sets Syst., {\bf 282} (2016), 1--20.\\
https://doi.org/10.1016/j.fss.2014.12.012
 \bibitem{Yao-Lu}W. Yao, L. X. Lu, {\it Fuzzy Galois connections on fuzzy posets}, Math. Logic Q., {\bf 55} (2009), 84--91.  https://doi.org/10.1002/malq.200710079
\bibitem{Yue-Yao-Ho}Y. L. Yue, W. Yao, W. K. Ho, {\it Applications of Scott-closed sets in convex structures}, Topol.  Appl., {\bf 314} (2022), 108093.\\
 https://doi.org/10.1016/j.topol.2022.108093
\bibitem{Yao-Zhou}W. Yao, C. J. Zhou, {\it Representation of sober convex spaces by join-semilattices}, J. Nonlinear Convex Anal., {\bf 21} (2020), 2715--2724.
 \bibitem{ZX-Zhang} Z. X. Zhang, {\it A general categorical reflection for various completions of closure spaces and $\mathcal{Q}$-ordered sets}, Fuzzy Sets Syst. {\bf 473} (2023), 1--18. https://doi.org/10.1016/j.fss.2023.108736
\bibitem{FanZhang}Q. Y. Zhang, L. Fan, {\it Continuity in quantitative domains},  Fuzzy Sets Syst., {\bf 154} (2005), 118--131. https://doi.org/10.1016/j.fss.2005.01.007
\bibitem{Zhao-Zhang} B. Zhao, Y. J. Zhang, Epireflections in the category of $T_0$ stratified $Q$-cotopological spaces, Soft Comput., {\bf 27} (2023), 2443--2451. \\
https://doi.org/10.1007/s00500-022-07809-y

\end{thebibliography}
\end{document}